\documentclass[preprint,12pt,1p]{elsarticle}
\makeatletter
\def\ps@pprintTitle{%
  \let\@oddhead\@empty
  \let\@evenhead\@empty
  \let\@oddfoot\@empty
  \let\@evenfoot\@oddfoot
}
\makeatother
\usepackage{graphics}
\usepackage{epsfig}
\usepackage{mathptmx}
\usepackage{amsmath}
\usepackage{amssymb}
\usepackage{hyperref}
\usepackage{fullpage}
\usepackage{amsthm}
\usepackage{geometry}
\usepackage{bigints}
\theoremstyle{plain}

\theoremstyle{definition}

\theoremstyle{example}

\theoremstyle{remark}

\newlength{\defbaselineskip}
\journal{Computer Aided Geometric Design}
\numberwithin{equation}{section}
\begin{document}
\begin{frontmatter}

\title{ Generalizations of Ramanujan's integral associated with infinite Fourier cosine transforms in terms of hypergeometric functions and its applications}
\author{M. I. Qureshi}
\author{*Showkat Ahmad Dar}
\ead{showkat34@gmail.com}
\address{E-Mail: miqureshi\_delhi@yahoo.co.in, showkat34@gmail.com \\Department of Applied Sciences and Humanities , \\Faculty of  Engineering and Technology , \\Jamia Millia Islamia (A Central University), New Delhi, 110025, India.}
\cortext[cor1]{Corresponding author}
\begin{abstract}
~~~~~In this paper, we obtain analytical solution of an unsolved integral $\textbf{R}_{C}(m,n)$ of Srinivasa Ramanujan [$\textit{Mess. Math}$., XLIV, 75-86, 1915], using hypergeometric approach, Mellin transforms, Infinite Fourier cosine transforms, Infinite series decomposition identity and some algebraic  properties of Pochhammer's symbol. Also we have given some generalizations of the Ramanujan's integral $\textbf{R}_{C}(m,n)$ in the form of integrals $\textbf{I}^{*}_{C}(\upsilon,b,c,\lambda,y), \textbf{J}_{C}(\upsilon,b,c,\lambda,y), \textbf{K}_{C}(\upsilon,b,c,\\\lambda,y),  \textbf{I}_{C}(\upsilon,b,\lambda,y)$ and solved it in terms of ordinary hypergeometric functions ${}_{2}F_{3}$, with suitable convergence conditions. Moreover as applications of Ramanujan's integral $\textbf{R}_{C}(m,n)$, the new nine infinite summation formulas associated with hypergeometric functions
${}_{0}F_{1}$, ${}_{1}F_{2}$ and ${}_{2}F_{3}$ are obtained.\\
\\
\textit{2010 AMS Classification: 33C20; 42A38; 33E20; 33C60.}
 \end{abstract}
\begin{keyword}
\small{ Generalized hypergeometric function; Infinite Fourier cosine transforms; Ramanujan's integrals; Fox-Wright psi hypergeometric function; Mellin transforms; Series decomposition identity; Bounded sequence.}
\end{keyword}
\end{frontmatter}
\section{Introduction and Preliminaries}
 In the literature of infinite Fourier cosine transforms  \cite{A15,B3,C1,C4,D,E2,M2,o,s3,s4,T1,W1}, the analytical solutions of  $\displaystyle{\int_{0}^{\infty}\frac{x^{\upsilon-1}\cos(xy)}{\{\exp(bx)\pm1\}}dx}$,
    are available in terms of Riemann's zeta function, the Psi function (Digamma function), hyperbolic function and Beta function. \\
The analytical solution of the following integral of Ramanujan \cite[p. 85, eq.(49) last line]{R1}:
\begin{equation}\label{b19}
\textbf{R}_{C}(m,n)= \int_{0}^{\infty}x^{m}\frac{\cos(\pi nx)}{\{-1+\exp{(2\pi\sqrt{x})}\}}dx,
\end{equation}
is not given for all positive rational values of $n$, and non-negative integral values of $m$.\\
For particular values of $m$ and $n$ in Ramanujan's integral $\textbf{R}_{C}(m,n)$, the following three integrals are given \cite[p.86, eq.(50)]{R1}:
\begin{equation}\label{s1}
\textbf{R}_{C}(1,1/2)=\int_{0}^{\infty}\frac{x\cos(\frac{\pi x}{2})}{\{-1+\exp{(2\pi\sqrt{x})}\}}dx=\frac{13-4\pi}{8\pi^{2}},~~~~~~~~~~~~~~~~
\end{equation}
\begin{equation}\label{s2}
\textbf{R}_{C}(1,2)=\int_{0}^{\infty}\frac{x\cos(2\pi x)}{\{-1+\exp{(2\pi\sqrt{x})}\}}dx=\frac{1}{64}\left(\frac{1}{2}-\frac{3}{\pi}+\frac{5}{\pi^{2}}\right),
\end{equation}
\begin{equation}\label{s3}
\textbf{R}_{C}(2,2)=\int_{0}^{\infty}\frac{x^{2}\cos(2\pi x)}{\{-1+\exp{(2\pi\sqrt{x})}\}}dx=\frac{1}{256}\left(1- \frac{5}{\pi}+\frac{5}{\pi^{2}}\right).
\end{equation}
The following theorem is proved by Ramanujan \cite[p.76-77, eqns(10 and $10^{'}$)]{R1}:\\
If
\begin{equation}\label{P1}
\textbf{R}_{C}(0,n)=\Phi (n)=\int_{0}^{\infty}\frac{\cos(\pi nx)}{\{-1+\exp{(2\pi\sqrt{x})}\}}dx,
\end{equation}
and
\begin{equation}\label{P2}
\Upsilon(n)=\frac{1}{2\pi n}+\int_{0}^{\infty}\frac{\sin(\pi nx)}{\{-1+\exp{(2\pi\sqrt{x})}\}}dx,
\end{equation}
then
\begin{equation}\label{P3}
\textbf{R}_{C}(0,n)=\Phi(n)=\frac{1}{n}\sqrt{\left(\frac{2}{n}\right)}~\Upsilon\left(\frac{1}{n}\right)-\Upsilon(n),
\end{equation}
and
\begin{equation}\label{P4}
\Upsilon(n)=\frac{1}{n}\sqrt{\left(\frac{2}{n}\right)}~\Phi\left(\frac{1}{n}\right)+\Phi(n),
\end{equation}
where $n$ is positive rational number.\\
For particular values of $n$, some values of Ramanujan's integral $\textbf{R}_{C}(0,n)=\Phi (n)$ \cite[p.85 (eq. 48)]{R1} are given below
\begin{equation}\label{P6}
\textbf{R}_{C}(0,1)=\Phi(1)=\int_{0}^{\infty}\frac{\cos(\pi x)}{\{-1+\exp{(2\pi\sqrt{x})}\}}dx=\frac{2-\sqrt{2}}{8},
\end{equation}
\begin{equation}\label{P7}
\textbf{R}_{C}(0,2)=\Phi(2)=\int_{0}^{\infty}\frac{\cos(2\pi x)}{\{-1+\exp{(2\pi\sqrt{x})}\}}dx=\frac{1}{16},
\end{equation}
\begin{equation}\label{P8}
\textbf{R}_{C}(0,4)=\Phi(4)=\int_{0}^{\infty}\frac{\cos(4\pi x)}{\{-1+\exp{(2\pi\sqrt{x})}\}}dx=\frac{3-\sqrt{2}}{32},
\end{equation}
\begin{equation}\label{P9}
\textbf{R}_{C}(0,6)=\Phi(6)=\int_{0}^{\infty}\frac{\cos(6\pi x)}{\{-1+\exp{(2\pi\sqrt{x})}\}}dx=\frac{13-4\sqrt{3}}{144},
\end{equation}
\begin{equation}\label{P10}
\textbf{R}_{C}(0,1/2)=\Phi\left(\frac{1}{2}\right)=\int_{0}^{\infty}\frac{\cos\left(\frac{\pi x}{2}\right)}{\{-1+\exp{(2\pi\sqrt{x})}\}}dx=\frac{1}{4\pi},
\end{equation}
\begin{equation}\label{P11}
\textbf{R}_{C}(0,2/5)=\Phi\left(\frac{2}{5}\right)=\int_{0}^{\infty}\frac{\cos\left(\frac{2\pi x}{5}\right)}{\{-1+\exp{(2\pi\sqrt{x})}\}}dx=\frac{8-3\sqrt{5}}{16}.
\end{equation}
\\
A natural generalization of Gauss hypergeometric series ${}_{2}F_{1}$ is the general hypergeometric series ${}_{p}F_{q}$  \cite[p.42, eq.(1)]{s1} and see also \cite{E1} with $p$ numerator parameters $\alpha_{1},  ... , \alpha_{p}$ and $q$ denominator parameters $\beta_{1}, ..., \beta_{q}$. It is defined by\\
\begin{equation}\label{h1}
{}_{p}F_{q}\left(\ \begin{array}{lll}\alpha_{1},...,\alpha_{p}~;~\\\beta_{1}, ..., \beta_{q}~;~\end{array} z\right)
=\sum_{n=0}^{\infty}\frac{(\alpha_{1})_{n}  ...  (\alpha_{p})_{n}}{(\beta_{1})_{n} ... (\beta_{q})_{n}}\frac{z^{n}}{n!}~,\newline
\end{equation}
where $\alpha_{i}\in\mathbb{C}~(i=1,...,p)$ and $\beta_{j}\in\mathbb{C}\setminus \mathbb{Z}_{0}^{-}~(j=1,...,q)~\left(\ \mathbb{Z}_{0}^{-}:=\{0,-1,-2,...\}\right)$ and \\
$\left(\ p,~q\in\mathbb{N}_{0}:=\mathbb{N}\cup\{0\}=\{0,1,2,...\}\right)$.
Also $(\lambda)_{\upsilon}~(\lambda, \upsilon\in\mathbb{C})$ denotes the Pochhammer's symbol (or the shifted factorial,
 since $(1)_{n}=n! )$ is defined, in general, by
\begin{equation}\label{L}
(\lambda)_{\upsilon}:=\frac{\Gamma(\lambda+\upsilon)}{\Gamma(\lambda)}=\begin{cases} 1 , \quad~~~~~~~~~~~~~ (\upsilon=0~;~\lambda\in\mathbb{C}\backslash\{0\}) \\ \lambda(\lambda+1)...(\lambda+n-1),  \quad (\upsilon=n\in\mathbb{N}~;~\lambda\in\mathbb{C}). \\
 \end{cases}
 \end{equation}
The hypergeometric ${}_{p}F_{q}$ series in the eq.(\ref{h1}) is convergent for $|z|<\infty$ if $p\leq q$, and for $|z|<1$ if $p=q+1$.\\
Furthermore, if we set
\begin{equation}\label{w}
\omega=\left(\ \sum_{j=1}^{q}\beta_{j}-\sum_{i=1}^{p}\alpha_{i}\right),\newline
\end{equation}
it is known that the ${}_{p}F_{q}$ series, with $p=q+1$, is\\
(i) absolutely convergent for $|z|=1$ if $Re(\omega)>0,$\\
(ii) conditionally convergent for $|z|=1, z\neq1$, if $-1<Re(\omega)\leq0$.\\
\\
Also binomial function is given by
\begin{equation}\label{b1}
(1-z)^{-\alpha}={}_{1}F_{0}\left(\begin{array}{lll}\alpha~;\\  \overline{~~~};\end{array} z \right)=\sum_{n=0}^{\infty}\frac{(\alpha)_{n}}{n!}z^{n}~,
\end{equation}
where $|z|<1,~\alpha\in\mathbb{C}\setminus \mathbb{Z}_{0}^{-}$,\\
The Fox-Wright psi function of one variable \cite{K5,K6,W2,W4} is given by
\begin{eqnarray}\label{F}
{}_{p}\Psi_{q}\left[\ \begin{array}{lll}(\alpha_{1}, A_{1}),...,(\alpha_{p}, A_{p});\\(\beta_{1}, B_{1}),...,(\beta_{q}, B_{q});~\end{array} z \right]
=\sum_{k=0}^{\infty}\frac{\Gamma(\alpha_{1}+kA_{1})...\Gamma(\alpha_{p}+kA_{p})}{\Gamma(\beta_{1}+kB_{1})...\Gamma(\beta_{q}+kB_{q})}\frac{z^{k}}{k!},
\end{eqnarray}
\begin{eqnarray}
~~~~~~~~~~~~~~~~~~~~~~~~~~~~~~~~~~~~=\frac{\Gamma(\alpha_{1})...\Gamma(\alpha_{p})}{\Gamma(\beta_{1})...\Gamma(\beta_{q})}\sum_{k=0}^{\infty}\frac{(\alpha_{1})_{kA_{1}}...(\alpha_{p})_{kA_{p}}}{(\beta_{1})_{kB_{1}}...(\beta_{q})_{kB_{q}}}\frac{z^{k}}{k!},\nonumber
\end{eqnarray}
\begin{eqnarray}\label{F3}
~~~~~~~~~~~~~~~~~~~~~~~~~~~~~~~~~~~~~~~~~~=\frac{\Gamma(\alpha_{1})...\Gamma(\alpha_{p})}{\Gamma(\beta_{1})...\Gamma(\beta_{q})}{}_{p}\Psi_{q}^{*}\left[\ \begin{array}{lll}(\alpha_{1}, A_{1}),...,(\alpha_{p}, A_{p});\\(\beta_{1}, B_{1}),...,(\beta_{q}, B_{q});~\end{array} z \right],
\end{eqnarray}
\begin{eqnarray}\label{F4}
~~~~~~~~~~~~~~~~~~~~~~~~~~~~~~~~~~~~~~~~~~~~~~~~~~~~~~~~~=\frac{1}{2\pi \rho}\bigint_{L}
\frac{\displaystyle \Gamma(\zeta)\prod_{i=1}^{p}\Gamma(\alpha_{i}-A_{i}\zeta)}{\displaystyle\prod_{j=1}^{q}\Gamma(\beta_{j}-B_{j}\zeta)}(-z)^{-\zeta}d\zeta~,~~~~~~~~~~~~~~
\end{eqnarray}
where $\rho=\sqrt{-1},~z\in\mathbb{C};~$ parameters $\alpha_{i},\beta_{j}\in\mathbb{C}$; coefficients $A_{i},B_{j}\in\mathbb{R}=(-\infty,+\infty)$ in case of series (\ref{F})  (or $A_{i},B_{j}\in\mathbb{R}_{+}=(0,+\infty)$ in case of contour integral (\ref{F4})), $A_{i}\neq0~(i=1,2,...,p),B_{j}\neq0~(j=1,2,...,q)$. In eq. (\ref{F}), the parameters $\alpha_{i},\beta_{j}$ and coefficients $A_{i},B_{j}$ are adjusted in such a way that the product of Gamma functions in numerator and denominator should be well defined.\\
Suppose:
\begin{equation}\label{W2}
\Delta^{*}=\left(\ \sum_{j=1}^{q}B_{j}-\sum_{i=1}^{p}A_{i}\right),
\end{equation}
\begin{equation}\label{W3}
\delta^{*}=\left(\prod_{i=1}^{p}|A_{i}|^{-A_{i}}\right)\left(\prod_{j=1}^{q}|B_{j}|^{B_{j}}\right),
\end{equation}
\begin{equation}\label{W4}
 \mu^{*}=\sum_{j=1}^{q}\beta_{j}-\sum_{i=1}^{p}\alpha_{i}+\left(\frac{p-q}{2}\right),
\end{equation}
and
\begin{equation}\label{W5}
\sigma^{*}=(1+A_{1}+...+A_{p})-(B_{1}+...+B_{q})=1-\Delta^{*}.
\end{equation}
Then we have the following convergence conditions of (\ref{F}) and (\ref{F4}):\\
\textbf{Case(1)}: When contour $(L)$ is a left loop beginning and ending at $-\infty$, then ${}_{p}\Psi_{q}[\cdot]$ given by (\ref{F})or (\ref{F4}) holds the following convergence conditions.\\
i)  When $\Delta^{*}>-1,$ $0<|z|<\infty$, $z\neq0$.\\
ii) When $\Delta^{*}=-1$, $0<|z|<\delta^{*}$. \\
iii) When $\Delta^{*}=-1$,$|z|=\delta^{*}$, and $Re(\mu^{*})>\frac{1}{2}$.\\
\textbf{Case(2)}: When contour $(L)$ is a right loop beginning and ending at $+\infty$, then ${}_{p}\Psi_{q}[\cdot]$ given by (\ref{F})or (\ref{F4}) holds the following convergence conditions.\\
iv) When $\Delta^{*}<-1$,  $0<|z|<\infty~,z\neq0$.\\ 
v) When $\Delta^{*}=-1$, $|z|>\delta^{*}$.\\ 
vi) When $\Delta^{*}=-1$,$|z|=\delta^{*}$, and $Re(\mu^{*})>\frac{1}{2}$.\\
\textbf{Case(3)}: When contour $(L)$ is starting from $\gamma-i\infty$ and ending at $\gamma+i\infty$ where $\gamma\in\mathbb{R}=(-\infty,+\infty)$, then ${}_{p}\Psi_{q}[\cdot]$ is also convergent under the following conditions.\\
vii) When $\sigma^{*}>0$, $|arg(-z)|<\frac{\pi}{2}\sigma^{*}$,  $0<|z|<\infty~,z\neq0$.\\
viii)  When $\sigma^{*}=0$,  $arg(-z)=0$ , $0<|z|<\infty,~z\neq0$ such that $-\gamma \Delta^{*}+Re(\mu^{*})>\frac{1}{2}+\gamma$.\\
ix) When $\gamma=0$, $\sigma^{*}=0, arg(-z)=0$, $0<|z|<\infty,~z\neq0$, such that $Re(\mu^{*})>\frac{1}{2}$.\\
The infinite Fourier cosine transform of $g(x)$ over the interval $(0,\infty)$ is defined by\\
\begin{equation}\label{f1}
F_{C}\{g(x);y\}=\int_{0}^{\infty}g(x)\cos(xy)dx=G_{C}(y),~~~~~(y>0),
\end{equation}
then $\displaystyle{g(x)=F_{C}^{-1}[G_{C}(y);x]=\frac{2}{\pi}\int_{0}^{\infty}G_{C}(y)\cos(xy)dy}$. \\
 Some authors have given the following definitions of infinite Fourier cosine transforms and its inverse.\\
$~~~~~~~~~~~~~~~~\mathcal{F}_{C}\{g(x);y\}=\displaystyle{\sqrt{\frac{2}{\pi}}\int_{0}^{\infty}g(x)\cos(xy)dx=\mathcal{G}_{C}(y),~~~~~(y>0)}$,\\
then $\displaystyle{g(x)=\mathcal{F}_{C}^{-1}[\mathcal{G}_{C}(y);x]=\sqrt{\frac{2}{\pi}}\int_{0}^{\infty}\mathcal{G}_{C}(y)\cos(xy)dy}$. \\
In this paper we shall apply the definition (\ref{f1}). If $b>0$ and $0<Re(s)<1$, then the Mellin-transform of $\cos(bx)$ or infinite Fourier Cosine transform of $x^{s-1}$\cite[p. 319, Entry 6.5(21), see also p.10, Entry 1.3(1)]{E2}, \cite[p.42, Entry(5.2)]{O1}, \cite[p.11,Entry(3.7)]{o} is given by
\begin{equation}\label{n1}
\int_{0}^{\infty}x^{s-1}\cos(bx)dx=\frac{\Gamma(s)\cos(\frac{\pi s}{2})}{b^{s}}.
\end{equation}
If $Re(\mu)>-1,$ $0<\xi<1$, $a>0$ and $y>0$, then we can prove the following integral by using Maclaurin's expansion of $\exp(-a x^{\xi})$ and term by term integrating with the help of the result (\ref{n1}):
\begin{equation}\label{n3}
\int_{0}^{\infty}x^{\mu}\exp(-ax^{\xi})\cos(xy)dx=-y^{-\mu-1}\sum_{\ell=0}^{\infty}\left(-\frac{a}{y^{\xi}}\right)^{\ell}\frac{1}{\ell!}
\Gamma(\mu+1+\xi\ell)~\sin\left\{\frac{\pi}{2}(\mu+\xi\ell)\right\}.
\end{equation}
The condition $Re(\mu)>-1$ stated in the integral (\ref{n3}) follows from the theory of analytic continuation \cite[p.15, Entry(3.55)]{o}, \cite[p.48, Entry(5.36)]{O1}.
We have also  verified the condition $Re(\mu)>-1$ , using Wolfram Mathematica.\\
An infinite series decomposition identity \cite[p.193,eq.(8)]{S5} and \cite{C2,L1,M1,S6,S7} is given by
\begin{equation}\label{K}
\sum_{\ell=0}^{\infty}\Omega(\ell)=\sum_{j=0}^{N-1}\left\{\sum_{\ell=0}^{\infty}\Omega(N\ell+j)\right\},
\end{equation}
where $N$ is an arbitrary positive integer.\\
Put $N=4$ in the above eq. (\ref{K}), we get
\begin{equation}\label{S1}
\sum_{\ell=0}^{\infty}\Omega(\ell)=\sum_{j=0}^{3}\left\{\sum_{\ell=0}^{\infty}\Omega(4\ell+j)\right\},~~~~~~~~~~~~~~~~~~~~~~~~~~~~~~~~~~~~~~~~~~~~~~~~~~~~~~~~~~~~~~~~
\end{equation}
\begin{equation}\label{S}
=\sum_{\ell=0}^{\infty}\Omega(4\ell)+\sum_{\ell=0}^{\infty}\Omega(4\ell+1)+\sum_{\ell=0}^{\infty}\Omega(4\ell+2)+
\sum_{\ell=0}^{\infty}\Omega(4\ell+3),
\end{equation}
provided that all involved infinite series are absolutely convergent.\\
For every positive integer $m$ \cite[p.22, eq.(26)]{s1}, we have\\
\begin{equation}\label{D}
  (\lambda)_{mn}=m^{mn}\prod_{j=1}^{m}\left(\frac{\lambda+j-1}{m}\right)_{n}~~~~~~~~;m\in\mathbb{N},~n\in\mathbb{N}_{0}\newline.
\end{equation}
From the above result (\ref{D}) with $\lambda=mz$, it can be proved that\\
\begin{equation}\label{M1}
  \Gamma(mz)=(2\pi)^{\frac{(1-m)}{2}}m^{mz-\frac{1}{2}}\prod_{j=1}^{m}\Gamma\left(z+\frac{j-1}{m}\right)\newline,~~~mz\in\mathbb{C}\backslash \mathbb{Z}_{0}^{-}
\end{equation}
~
The equation (\ref{M1}) is known as Gauss-Legendre multiplication theorem for Gamma function.\\
Elementary trigonometric functions \cite[p.44, eq.(9) and eq.(10)]{s1} are given by
\begin{equation}\label{E}
\cos~z={}_{0}F_{1}\left(\ \begin{array}{lll}\overline{~~~~~};~\\ ~~~\frac{1}{2};~\end{array} \frac{-z^{2}}{4}\right),
\end{equation}
\begin{equation}\label{E1}
\sin~z=z~{}_{0}F_{1}\left(\ \begin{array}{lll}\overline{~~~~~};~\\ ~~~\frac{3}{2};~\end{array} \frac{-z^{2}}{4}\right).
\end{equation}
Lommel function \cite[p.44, eq.(13)]{s1} is given by
\begin{equation}\label{E2}
s_{\mu,\upsilon}(z)=\frac{z^{\mu+1}}{(\mu-\upsilon+1)(\mu+\upsilon+1)}{}_{1}F_{2}\left(\ \begin{array}{lll}~~~~~~~~~~~~~~~~~~~~~~1;~\\ \frac{\mu-\upsilon+3}{2},\frac{\mu+\upsilon+3}{2};~\end{array} \frac{-z^{2}}{4}\right),
\end{equation}
where $\mu\pm\upsilon\in\mathbb{C}\setminus \{-1,-3,-5,-7,...\}$.\\
\\
~~~~In the available literature \cite{A2,A3,A4,A11,A12,A13,A14,B10,B11,B12,B13,B14,P2,P3,P4,P5,P6,R1,R4,R0,R5,R3,R2} on Ramanujan's Mathematics, the analytical solution of Ramanujan's integral $\textbf{R}_{C}(m,n)$ is not given. Therefore, the main aim of this paper is to obtain the analytical solution of Ramanujan's integral in terms of ordinary hypergeometric functions. Also, our work on Ramanujan's Mathematics is motivated by the work given in references \cite{M,Q12,Q9,Q4,Q2,s2}.
Here in this paper, we generalized Ramanujan's integral $\textbf{R}_{C}(m,n)$  in the following forms: \\
 (i)~ $\textbf{I}^{*}_{C}(\upsilon,b,c,\lambda,y)=\displaystyle{\sum_{k=0}^{\infty}}\bigg[\frac{\Theta(k)}{k!}
\int_{0}^{\infty}~x^{\upsilon-1}e^{-(\lambda b +c k)\sqrt{x}}\cos (xy)dx\bigg]$,\\
(ii)~ $\textbf{J}_{C}(\upsilon,b,c,\lambda,y)=\displaystyle{\int_{0}^{\infty}}x^{\upsilon-1}e^{-b\lambda\sqrt{x}}
{}_{r}\Psi_{s}\left[\ \begin{array}{lll}(\alpha_{1}, A_{1}),...,(\alpha_{r}, A_{r});\\(\beta_{1}, B_{1}),...,(\beta_{s}, B_{s});~\end{array} e^{-c\sqrt{x}}\right]\cos (xy)dx$,\\
(iii)~$\textbf{K}_{C}(\upsilon,b,c,\lambda,y)=\displaystyle{\int_{0}^{\infty}}x^{\upsilon-1}e^{-b\lambda\sqrt{x}}{}_{r}F_{s}\left(\ \begin{array}{lll}\alpha_{1},...,\alpha_{r};\beta_{1},...,\beta_{s};~\end{array} e^{-c\sqrt{x}}\right)\cos (xy)dx$,\\
(iv)~ $\textbf{I}_{C}(\upsilon,b,\lambda,y)=
\displaystyle{\int_{0}^{\infty}}x^{\upsilon-1}
\{-1+exp(b\sqrt{x})\}^{-\lambda}cos(xy)dx$, \\where $\{\Theta(k)\}_{k=0}^{\infty}$ is a bounded sequence and obtained the analytical solution. Moreover, we also show how the main general theorem given below, is applicable for obtaining new interesting results  by suitable adjustment in parameters and variables (see in the sections 3,4,5,6).
\section{Main general theorem on infinite Fourier cosine transform}
 Suppose $\{\Theta(k)\}_{k=0}^{\infty}$ is a bounded sequence of arbitrary real and complex numbers, and $Re(\upsilon)>0, c>0,y>0;\lambda>0,b>0 (or~\lambda<0,b<0)$ then
\begin{equation}\label{L1}
\textbf{I}^{*}_{C}(\upsilon,b,c,\lambda,y)=\sum_{k=0}^{\infty}\left[\frac{\Theta(k)}{k!}
\int_{0}^{\infty}~x^{\upsilon-1}e^{-(\lambda b +c k)\sqrt{x}} \cos (xy)dx\right],~~~~~~~~~~~~~~~~~~~~~~~~~~~~~~~~~~~~
\end{equation}
\begin{equation}\label{L3}
=y^{-\upsilon}\sum_{k=0}^{\infty}\bigg[\frac{\Theta(k)}{k!}\sum_{\ell=0}^{\infty}\frac{(-1)^{\ell}(\lambda b+c k)^{\ell}\Gamma\left(\upsilon+\frac{\ell}{2}\right)}{y^{\frac{\ell}{2}}~\ell!}\cos\left(\frac{\upsilon \pi}{2}+\frac{\ell\pi}{4}\right)\bigg],
\end{equation}
Now replacing $\ell$ by $4\ell+j$, after simplification we get
\begin{multline}\label{L4}
\textbf{I}^{*}_{C}(\upsilon,b,c,\lambda,y)=y^{-\upsilon}\sum_{k=0}^{\infty}\bigg[\frac{\Theta(k)}{k!}\sum_{j=0}^{3}\frac{(-1)^{j}(\lambda b+c k)^{j}~\Gamma\left(\upsilon+\frac{j}{2}\right)}{y^{\frac{j}{2}}~j!}\cos\left(\frac{\upsilon \pi}{2}+\frac{j\pi}{4}\right)\times\\
\times {}_{2}F_{3}\left(\
\begin{array}{lll} \Delta\left(2;\frac{2\upsilon+j}{2}\right);\\ \Delta^{*}\left(4;1+j\right);\end{array} \frac{-1}{64y^{2}}\left\{\frac{(\lambda b)(\frac{\lambda b+c}{c})_{k}}{(\frac{\lambda b}{c})_{k}}\right\}^{4}\right)\bigg],
\end{multline}
\begin{multline}\label{L5}
~~~~~~~~~~~~~~~~~~~~~~=y^{-\upsilon}\sum_{k=0}^{\infty}\bigg[\frac{\Theta(k)}{k!}\sum_{j=0}^{3}\frac{(-1)^{j}\Gamma\left(\upsilon+\frac{j}{2}\right)}{j!}\cos\left(\frac{\upsilon \pi}{2}+\frac{j\pi}{4}\right)\left(\frac{\lambda b}{\sqrt{y}}\right)^{j}\times\\\times\left\{\frac{(\frac{\lambda b+c}{c})_{k}}{(\frac{\lambda b}{c})_{k}}\right\}^{j}
 {}_{2}F_{3}\left(\
\begin{array}{lll} \Delta\left(2;\frac{2\upsilon+j}{2}\right);\\ \Delta^{*}\left(4;1+j\right);\end{array} \frac{-1}{64y^{2}}\left\{\frac{(\lambda b)(\frac{\lambda b+c}{c})_{k}}{(\frac{\lambda b}{c})_{k}}\right\}^{4}\right)\bigg],
\end{multline}
\begin{eqnarray}\label{L2}
~~~~~~~~~~~~~~~~~~~~~~~~=\frac{\Gamma(\upsilon)\cos\left(\frac{\upsilon \pi}{2}\right)}{y^{\upsilon}}\sum_{k=0}^{\infty}\bigg[\frac{\Theta(k)}{k!}{}_{2}F_{3}\left(\
\begin{array}{lll} \frac{\upsilon}{2},\frac{\upsilon+1}{2}~;\\ \frac{1}{4},~\frac{1}{2},~\frac{3}{4};\end{array} \frac{-1}{64y^{2}}\left\{\frac{(\lambda b)(\frac{\lambda b+c}{c})_{k}}{(\frac{\lambda b}{c})_{k}}\right\}^{4}\right)\bigg]\nonumber-
\end{eqnarray}
\begin{multline}
-\frac{(\lambda b)\Gamma(\upsilon+\frac{1}{2})\cos\left(\frac{\upsilon \pi}{2}+\frac{\pi}{4}\right)}{y^{\upsilon+\frac{1}{2}}}\sum_{k=0}^{\infty}\bigg[\frac{\Theta(k)}{k!}\left\{\frac{(\frac{\lambda b+c}{c})_{k}}{(\frac{\lambda b}{c})_{k}}\right\}\times\\\times
{}_{2}F_{3}\left(\
\begin{array}{lll} \frac{2\upsilon+1}{4},\frac{2\upsilon+3}{4};\\ \frac{1}{2},~\frac{3}{4},~\frac{5}{4}~~~~~;\end{array} \frac{-1}{64y^{2}}\left\{\frac{(\lambda b)(\frac{\lambda b+c}{c})_{k}}{(\frac{\lambda b}{c})_{k}}\right\}^{4}\right)\bigg]\nonumber-
\end{multline}
\begin{multline}
-\frac{(\lambda b)^{2}\Gamma(\upsilon+1)\sin\left(\frac{\upsilon \pi}{2}\right)}{2y^{\upsilon+1}}\sum_{k=0}^{\infty}\bigg[\frac{\Theta(k)}{k!}\left\{\frac{(\frac{\lambda b+c}{c})_{k}}{(\frac{\lambda b}{c})_{k}}\right\}^{2}\times\\\times
{}_{2}F_{3}\left(\
\begin{array}{lll} \frac{\upsilon+1}{2},\frac{\upsilon+2}{2};\\ \frac{3}{4},~\frac{5}{4},~\frac{3}{2}~~;\end{array} \frac{-1}{64y^{2}}\left\{\frac{(\lambda b)(\frac{\lambda b+c}{c})_{k}}{(\frac{\lambda b}{c})_{k}}\right\}^{4}\right)\bigg]\nonumber+
\end{multline}
\begin{multline}\label{K3}
+\frac{(\lambda b)^{3}\Gamma(\upsilon+\frac{3}{2})\sin\left(\frac{\upsilon \pi}{2}+\frac{\pi}{4}\right)}{6y^{\upsilon+\frac{3}{2}}}\sum_{k=0}^{\infty}\bigg[\frac{\Theta(k)}{k!}\left\{\frac{(\frac{\lambda b+c}{c})_{k}}{(\frac{\lambda b}{c})_{k}}\right\}^{3}
\times\\\times{}_{2}F_{3}\left(\
\begin{array}{lll} \frac{2\upsilon+3}{4},\frac{2\upsilon+5}{4};\\ \frac{5}{4},~\frac{3}{2},~\frac{7}{4}~~~~~;\end{array} \frac{-1}{64y^{2}}\left\{\frac{(\lambda b)(\frac{\lambda b+c}{c})_{k}}{(\frac{\lambda b}{c})_{k}}\right\}^{4}\right)\bigg].
\end{multline}
Our result (\ref{L4})or (\ref{L5}) or (\ref{K3}) is convergent in view of the convergence condition of ${}_{p}F_{q}(\cdot)$ series, when $p\leq q$, and $ \forall ~|z|<\infty$.\\
\textbf{Proof}: The result (\ref{L3}) is obtained by the application of the integral (\ref{n3}) [with substitutions $\mu=\upsilon-1,a=\lambda b+ck,\xi=\frac{1}{2}$] in the R.H.S. of eq.(\ref{L1}).
The results (\ref{L4}), (\ref{L5}) and (\ref{L2}) are obtained by using the infinite series decomposition formulas (\ref{S1}),(\ref{S}), Pochhammer's identity (\ref{D}) and other algebraic properties of Pochhammer's symbols.
\section{Infinite Fourier cosine transforms of $x^{\upsilon-1}e^{-b\lambda\sqrt{x}}{}_{r}\Psi_{s}[\cdot]$ and  $x^{\upsilon-1}e^{-b\lambda\sqrt{x}}{}_{r}F_{s}(\cdot)$}
If we put $\Theta(k)=\displaystyle{\frac{\Gamma(\alpha_{1}+kA_{1})...\Gamma(\alpha_{r}+kA_{r})} {\Gamma(\beta_{1}+kB_{1})...\Gamma(\beta_{s}+kB_{s})}}$, $ (k=0,1,2,3,...)$
 in the equations (\ref{L1}) and (\ref{L4}), then after simplification we get (\ref{G}) and (\ref{G2})
\begin{equation}\label{G}
\textbf{J}_{C}(\upsilon,b,c,\lambda,y)=\int_{0}^{\infty}x^{\upsilon-1}e^{-b\lambda\sqrt{x}}{}_{r}\Psi_{s}\left[\ \begin{array}{lll}(\alpha_{1}, A_{1}),...,(\alpha_{r}, A_{r});\\(\beta_{1}, B_{1}),...,(\beta_{s}, B_{s});~\end{array} e^{-c\sqrt{x}}\right]\cos (xy)dx,
\end{equation}
\begin{multline}\label{G2}
~~~~~~~~~~~~~~~~~~~~~~~=y^{-\upsilon}\sum_{k=0}^{\infty}\bigg[\frac{\Gamma(\alpha_{1}+kA_{1})...\Gamma(\alpha_{r}+kA_{r})} {\Gamma(\beta_{1}+kB_{1})...\Gamma(\beta_{s}+kB_{s})k!}
\sum_{j=0}^{3}\frac{(-1)^{j}(\lambda b+c k)^{j}~\Gamma\left(\upsilon+\frac{j}{2}\right)}{y^{\frac{j}{2}}~j!}\times\\\times \cos\left(\frac{\upsilon \pi}{2}+\frac{j\pi}{4}\right)
 {}_{2}F_{3}\left(\
\begin{array}{lll} \Delta\left(2;\frac{2\upsilon+j}{2}\right);\\ \Delta^{*}\left(4;1+j\right);\end{array} \frac{-1}{64y^{2}}\left\{\frac{(\lambda b)(\frac{\lambda b+c}{c})_{k}}{(\frac{\lambda b}{c})_{k}}\right\}^{4}\right)\bigg],
\end{multline}
where  $Re(\upsilon)>0;~c>0,y>0;\lambda>0,b>0 ~(or~\lambda<0,b<0)$;~parameters $\alpha_{i},\beta_{j}\in\mathbb{C}$;~coefficients  $A_{i},B_{j}\in\mathbb{R}=(-\infty,+\infty);A_{i}\neq0~(i=1,2,...,r),B_{j}\neq0~(j=1,2,...,s)$ and ${}_{r}\Psi_{s}[\cdot]$ is Fox-Wright psi function of one variable subject to suitable convergence conditions derived from convergence conditions discussed in \textbf{case(1)} or \textbf{case(2)} or \textbf{case(3)} of the function ${}_{p}\Psi_{q}[\cdot]$ given by (\ref{F}),(\ref{F3}) and (\ref{F4}).\\
$~~~~~~~~~~~~~$When $N$ is positive integer then $\Delta(N;\lambda)$ denotes the array of $N$ parameters given by $\frac{\lambda}{N},\frac{\lambda+1}{N},...,\frac{\lambda+N-1}{N}$. When $N$ and $j$ are independent variables then the notation $\Delta(N;j+1)$ denotes the set of $N$ parameters given by $\frac{j+1}{N},\frac{j+2}{N},...,\frac{j+N}{N}$. When $j$ is dependent variable that is $j=0,1,2,3,...,N-1$, then the asterisk in $\Delta^{*}(N;j+1)$ represents the fact that the (denominator) parameters $\frac{N}{N}$ is always omitted (due to the need of factorial in denominator in the power series form of hypergeometric function) so that the set $\Delta^{*}(N;j+1)$ obviously contains only $(N-1)$ parameters \cite[Chap.3, p.214]{s1}.\\
\\
\textbf{Remark}: When $A_{1}=...=A_{r}=B_{1}=...=B_{s}=1$ in (\ref{G}), (\ref{G2}) then we get
\begin{equation}\label{G3}
\textbf{K}_{C}(\upsilon,b,c,\lambda,y)=\int_{0}^{\infty}x^{\upsilon-1}e^{-b\lambda\sqrt{x}}{}_{r}F_{s}\left(\ \begin{array}{lll}\alpha_{1},...,\alpha_{r};\\\beta_{1},...,\beta_{s};~\end{array} e^{-c\sqrt{x}}\right)\cos (xy)dx,
\end{equation}
\begin{multline}\label{G4}
~~~~~~~~~~~~~~~~~~=y^{-\upsilon}\sum_{k=0}^{\infty}\bigg[\frac{(\alpha_{1})_{k}...(\alpha_{r})_{k}} {(\beta_{1})_{k}...(\beta_{s})_{k}~k!}
\sum_{j=0}^{3}\frac{(-1)^{j}(\lambda b+c k)^{j}~\Gamma\left(\upsilon+\frac{j}{2}\right)}{y^{\frac{j}{2}}~j!} \cos\left(\frac{\upsilon \pi}{2}+\frac{j\pi}{4}\right)\times\\\times
 {}_{2}F_{3}\left(\
\begin{array}{lll} \Delta\left(2;\frac{2\upsilon+j}{2}\right);\\ \Delta^{*}\left(4;1+j\right);\end{array} \frac{-1}{64y^{2}}\left\{\frac{(\lambda b)(\frac{\lambda b+c}{c})_{k}}{(\frac{\lambda b}{c})_{k}}\right\}^{4}\right)\bigg],
\end{multline}
where  $Re(\upsilon)>0;~c>0,y>0;\lambda>0,b>0 ~(or~\lambda<0,b<0)$;~parameters $\alpha_{i},\beta_{j}\in\mathbb{C} ~(i=1,2,...,r),(j=1,2,...,s)$ and $r\leq s+1$.
\section{Infinite Fourier cosine transform of $x^{\upsilon-1}\{-1+\exp{(b\sqrt{x})}\}^{-\lambda}$}
The following generalization $\textbf{I}_{C}(\upsilon,b,\lambda,y)$ of the Ramanujan's integral $\textbf{R}_{C}(m,n)$ in terms of ordinary hypergeometric functions ${}_{2}F_{3}$ holds true:
\begin{equation}\label{D1}
\textbf{I}_{C}(\upsilon,b,\lambda,y)=\int_{0}^{\infty}x^{\upsilon-1}\frac{\cos(xy)}{\{-1+\exp{(b\sqrt{x})}\}^{\lambda}}dx,~~~~~~~~~~~~~~~~~~~~~~~~~~~~~~~~~~~~~~~~~~~~~~~~~~~~~~~~~
\end{equation}
\begin{equation}\label{D3}
=y^{-\upsilon}\sum_{k=0}^{\infty}\bigg[~\frac{(\lambda)_{k}}{k!}\sum_{\ell=0}^{\infty}\frac{(-1)^{\ell}(\lambda b+b k)^{\ell}\Gamma\left(\upsilon+\frac{\ell}{2}\right)}{y^{\frac{\ell}{2}}~\ell!}\cos\left(\frac{\upsilon \pi}{2}+\frac{\ell\pi}{4}\right)\bigg],~~~~~~~~~
\end{equation}
\begin{multline}\label{D4}
~~~~~~~~~~~~~~=y^{-\upsilon}\sum_{k=0}^{\infty}\bigg[\frac{(\lambda)_{k}}{k!}\sum_{j=0}^{3}\frac{(-1)^{j}(\lambda b+b k)^{j}~\Gamma\left(\upsilon+\frac{j}{2}\right)}{y^{\frac{j}{2}}~j!}\cos\left(\frac{\upsilon \pi}{2}+\frac{j\pi}{4}\right)\times\\
\times {}_{2}F_{3}\left(\
\begin{array}{lll} \Delta\left(2;\frac{2\upsilon+j}{2}\right);\\ \Delta^{*}\left(4;1+j\right);\end{array} \frac{-1}{64y^{2}}\left\{\frac{(\lambda b)(\lambda+1)_{k}}{(\lambda)_{k}}\right\}^{4}\right)\bigg],
\end{multline}
\begin{eqnarray}\label{D2}
~~~~~~=\frac{\Gamma(\upsilon)\cos\left(\frac{\upsilon \pi}{2}\right)}{y^{\upsilon}}\sum_{k=0}^{\infty}\bigg[\frac{(\lambda)_{k}}{k!}{}_{2}F_{3}\left(\
\begin{array}{lll} \frac{\upsilon}{2},\frac{\upsilon+1}{2}~;\\ \frac{1}{4},~\frac{1}{2},~\frac{3}{4};\end{array} \frac{-1}{64y^{2}}\left\{\frac{(\lambda b)(\lambda+1)_{k}}{(\lambda)_{k}}\right\}^{4}\right)\bigg]\nonumber-
\end{eqnarray}
\begin{eqnarray}
-\frac{(\lambda b)\Gamma(\upsilon+\frac{1}{2})\cos\left(\frac{\upsilon \pi}{2}+\frac{\pi}{4}\right)}{y^{\upsilon+\frac{1}{2}}}\sum_{k=0}^{\infty}\bigg[\frac{(\lambda+1)_{k}}{k!}{}_{2}F_{3}\left(\
\begin{array}{lll} \frac{2\upsilon+1}{4},\frac{2\upsilon+3}{4};\\ \frac{1}{2},~\frac{3}{4},~\frac{5}{4}~~~~~;\end{array} \frac{-1}{64y^{2}}\left\{\frac{(\lambda b)(\lambda+1)_{k}}{(\lambda)_{k}}\right\}^{4}\right)\bigg]\nonumber-\\
-\frac{(\lambda b)^{2}\Gamma(\upsilon+1)\sin\left(\frac{\upsilon \pi}{2}\right)}{2y^{\upsilon+1}}\sum_{k=0}^{\infty}\bigg[\frac{\left\{(\lambda+1)_{k}\right\}^{2}}{(\lambda)_{k}~k!}{}_{2}F_{3}\left(\
\begin{array}{lll} \frac{\upsilon+1}{2},\frac{\upsilon+2}{2};\\ \frac{3}{4},~\frac{5}{4},~\frac{3}{2}~~;\end{array} \frac{-1}{64y^{2}}\left\{\frac{(\lambda b)(\lambda+1)_{k}}{(\lambda)_{k}}\right\}^{4}\right)\bigg]\nonumber+\\
+\frac{(\lambda b)^{3}\Gamma(\upsilon+\frac{3}{2})\sin\left(\frac{\upsilon \pi}{2}+\frac{\pi}{4}\right)}{6y^{\upsilon+\frac{3}{2}}}\sum_{k=0}^{\infty}\bigg[\frac{\left\{(\lambda+1)_{k}\right\}^{3}}{k!\left\{(\lambda)_{k}\right\}^{2}}{}_{2}F_{3}\left(\
\begin{array}{lll} \frac{2\upsilon+3}{4},\frac{2\upsilon+5}{4};\\ \frac{5}{4},~\frac{3}{2},~\frac{7}{4}~~~~~;\end{array} \frac{-1}{64y^{2}}\left\{\frac{(\lambda b)(\lambda+1)_{k}}{(\lambda)_{k}}\right\}^{4}\right)\bigg], \nonumber\\
\end{eqnarray}
where $Re(\upsilon)>0;~y>0;\lambda>0,b>0$.\\
\\
\textbf{Proof}: In eq.(\ref{L1}), put $\Theta(k)=(\lambda)_{k}$ and $c=b$, we obtain
\begin{equation}\label{K1}
\textbf{I}_{C}(\upsilon,b,\lambda,y)=\int_{0}^{\infty}~x^{\upsilon-1}e^{-(\lambda b)\sqrt{x}}\left\{\sum_{k=0}^{\infty}\frac{(\lambda)_{k}}{k!}e^{-(bk)\sqrt{x}}\right\} \cos (xy)dx.
\end{equation}
Using binomial expansion (\ref{b1}) in the above eq. (\ref{K1}), after simplification, we get the equation (\ref{D1}).
The equations  (\ref{D3}), (\ref{D4}) and (\ref{D2}) are obtained from (\ref{L3}), (\ref{L4}) and (\ref{L2}) by putting $\Theta(k)=(\lambda)_{k}$ and $c=b$.
\section{Ramanujan's integral $\textbf{R}_{C}(m,n)$}
The analytical solution of the integral $\textbf{R}_{C}(m,n)$ is given by
\begin{equation}\label{K2}
\textbf{R}_{C}(m,n)=\int_{0}^{\infty}x^{m}\frac{\cos(\pi nx)}{\{-1+\exp{(2\pi\sqrt{x})}\}}dx,~~~~~~~~~~~~~~~~~~~~~~~~~~~~~~~~~~~~~~~~~~~~~~~~~~~~~~~~~~~~~~~~~~~~~~~~~~
\end{equation}
\begin{eqnarray}\label{N}
=-(n\pi)^{-m-1}\sum_{k=0}^{\infty}\bigg[\sum_{\ell=0}^{\infty}\frac{1}{\ell!}\left\{\frac{-(2\pi+2\pi k)}{\sqrt{n\pi}}\right\}^{\ell}\Gamma\left(m+1+\frac{\ell}{2}\right)\sin\left(\frac{m \pi}{2}+\frac{\ell\pi}{4}\right)\bigg],
 \end{eqnarray}
\begin{multline}\label{K4}
~~~=-(n\pi)^{-m-1}\sum_{k=0}^{\infty}\bigg[\sum_{j=0}^{3}\frac{1}{j!}\left\{\frac{-(2\pi+2\pi k)}{\sqrt{n\pi}}\right\}^{j}\Gamma\left(m+1+\frac{j}{2}\right)\sin\left(\frac{m \pi}{2}+\frac{j\pi}{4}\right)\times\\
\times {}_{2}F_{3}\left(\
\begin{array}{lll} \Delta\left(2;\frac{2m+j+2}{2}\right);\\ \Delta^{*}\left(4;1+j\right);\end{array} \frac{-\pi^{2}}{4n^{2}}\left\{\frac{(2)_{k}}{(1)_{k}}\right\}^{4}\right)\bigg],~~~~~~
\end{multline}
\begin{eqnarray}\label{F2}
=-\frac{m!~\sin\left(\frac{m \pi}{2}\right)}{(n\pi)^{m+1}}\sum_{k=0}^{\infty}\bigg[{}_{2}F_{3}\left(\
\begin{array}{lll} \frac{m+1}{2},\frac{m+2}{2};\\ \frac{1}{4},~\frac{1}{2},~\frac{3}{4}~~~;\end{array} -\frac{\pi^{2}}{4n^{2}}\left\{\frac{(2)_{k}}{(1)_{k}}\right\}^{4}\right)\bigg]\nonumber+\\
+\frac{\left(\frac{3}{2}\right)_{m}\sin\left(\frac{m \pi}{2}+\frac{\pi}{4}\right)}{(\pi)^{m}(n)^{m+\frac{3}{2}}}\sum_{k=0}^{\infty}\bigg[\left\{\frac{(2)_{k}}{(1)_{k}}\right\}{}_{2}F_{3}\left(\
\begin{array}{lll} \frac{2m+3}{4},\frac{2m+5}{4};\\ \frac{1}{2},~\frac{3}{4},~\frac{5}{4}~~~~~~;\end{array} \frac{-\pi^{2}}{4n^{2}}\left\{\frac{(2)_{k}}{(1)_{k}}\right\}^{4}\right)\bigg]\nonumber-\\
-\frac{(2)(m+1)!\cos\left(\frac{m \pi}{2}\right)}{(\pi)^{m}(n)^{m+2}}\sum_{k=0}^{\infty}\bigg[\left\{\frac{(2)_{k}}{(1)_{k}}\right\}^{2}{}_{2}F_{3}\left(\
\begin{array}{lll} \frac{m+2}{2},\frac{m+3}{2};\\ \frac{3}{4},~\frac{5}{4},~\frac{3}{2}~~~;\end{array} \frac{-\pi^{2}}{4n^{2}}\left\{\frac{(2)_{k}}{(1)_{k}}\right\}^{4}\right)\biggr]\nonumber+\\
+\frac{\left(\frac{5}{2}\right)_{m}\cos\left(\frac{m \pi}{2}+\frac{\pi}{4}\right)}{(\pi)^{m-1}(n)^{m+\frac{5}{2}}}\sum_{k=0}^{\infty}\bigg[\left\{\frac{(2)_{k}}{(1)_{k}}\right\}^{3}{}_{2}F_{3}\left(\
\begin{array}{lll} \frac{2m+5}{4},\frac{2m+7}{4};\\ \frac{5}{4},~\frac{3}{2},~\frac{7}{4}~~~~~~;\end{array} \frac{-\pi^{2}}{4n^{2}}\left\{\frac{(2)_{k}}{(1)_{k}}\right\}^{4}\right)\bigg],
\end{eqnarray}
where  $m$ is a non-negative integer and $n$ is positive rational number.\\
\\
\textbf{Proof}: The results (\ref{K2}), (\ref{N}), (\ref{K4}) and (\ref{F2}) are obtained from (\ref{D1}), (\ref{D3}), (\ref{D4}) and (\ref{D2}) by putting $\upsilon=m+1,~b=2\pi$,~$\lambda=1$ and $y=n\pi$.
\section{Applications of Ramanujan's integrals}
In this section we have established the following nine infinite new summation
formulas associated with hypergeometric series ${}_{0}F_{1}$, ${}_{1}F_{2}$ and ${}_{2}F_{3}$:
\begin{eqnarray}\label{N1}
\sum_{k=0}^{\infty}\bigg[{}_{2}F_{3}\left(\ \begin{array}{lll}1,\frac{3}{2}~~~~;\\\frac{1}{4},\frac{1}{2},\frac{3}{4};\end{array} -\pi^{2}\left\{\frac{(2)_{k}}{(1)_{k}}\right\}^{4}\right)\bigg]
-\frac{3\pi}{2}\sum_{k=0}^{\infty}\bigg[\left\{\frac{(2)_{k}}{(1)_{k}}\right\}~{}_{1}F_{2}\left(\ \begin{array}{lll}\frac{7}{4}~~~~;\\ \frac{1}{2},\frac{3}{4};\end{array} -\pi^{2}\left\{\frac{(2)_{k}}{(1)_{k}}\right\}^{4}\right)\bigg]\nonumber+\\
+5\pi^{2}\sum_{k=0}^{\infty}\bigg[\left\{\frac{(2)_{k}}{(1)_{k}}\right\}^{3}{}_{1}F_{2}\left(\ \begin{array}{lll}\frac{9}{4}~~~~;\\
\frac{5}{4},\frac{3}{2};\end{array} -\pi^{2}\left\{\frac{(2)_{k}}{(1)_{k}}\right\}^{4} \right)\bigg]
=\frac{1}{32}\left(4\pi-13\right),~~~~~~~~~~~~
\end{eqnarray}
\begin{eqnarray}\label{N2}
\sum_{k=0}^{\infty}\bigg[{}_{2}F_{3}\left(\
\begin{array}{lll} 1,\frac{3}{2}~~~~;\\ \frac{1}{4},\frac{1}{2},\frac{3}{4};\end{array} \frac{-\pi^{2}}{16}\left\{\frac{(2)_{k}}{(1)_{k}}\right\}^{4}\right)\bigg]
-\frac{3\pi}{4}\sum_{k=0}^{\infty}\bigg[\left\{\frac{(2)_{k}}{(1)_{k}}\right\}{}_{1}F_{2}\left(\ \begin{array}{lll}\frac{7}{4}~~~~;\\ \frac{1}{2},\frac{3}{4};\end{array} \frac{-\pi^{2}}{16}\left\{\frac{(2)_{k}}{(1)_{k}}\right\}^{4}\right)\bigg]\nonumber+\\
+\frac{5\pi^{2}}{8}\sum_{k=0}^{\infty}\bigg[\left\{\frac{(2)_{k}}{(1)_{k}}\right\}^{3}{}_{1}F_{2}\left(\ \begin{array}{lll}\frac{9}{4}~~~~;\\
\frac{5}{4},\frac{3}{2};~\end{array} \frac{-\pi^{2}}{16}\left\{\frac{(2)_{k}}{(1)_{k}}\right\}^{4}\right)\bigg]
=\frac{\pi^{2}}{16}\left(\frac{3}{\pi}-\frac{1}{2}-\frac{5}{\pi^{2}}\right),~~~~~~~~
\end{eqnarray}
\begin{align}\label{N3}
\sum_{k=0}^{\infty}\bigg[\left\{\frac{(2)_{k}}{(1)_{k}}\right\}{}_{2}F_{3}\left(\ \begin{array}{lll}\frac{7}{4},\frac{9}{4}~~~~;\\\frac{1}{2},\frac{3}{4},\frac{5}{4};~\end{array} \frac{-\pi^{2}}{16}\left\{\frac{(2)_{k}}{(1)_{k}}\right\}^{4}\right)\bigg]
-\frac{16}{5}\sum_{k=0}^{\infty}\bigg[\left\{\frac{(2)_{k}}{(1)_{k}}\right\}^{2}{}_{2}F_{3}\left(\ \begin{array}{lll}2,\frac{5}{2}~~~~;\\\frac{3}{4},\frac{5}{4},\frac{3}{2};\end{array} \frac{-\pi^{2}}{16}\left\{\frac{(2)_{k}}{(1)_{k}}\right\}^{4}\right)\bigg]\nonumber+\\
+\frac{7\pi}{6}\sum_{k=0}^{\infty}\bigg[\left\{\frac{(2)_{k}}{(1)_{k}}\right\}^{3}{}_{2}F_{3}\left(\ \begin{array}{lll}\frac{9}{4},\frac{11}{4}~~;\\\frac{5}{4},\frac{3}{2},\frac{7}{4};\end{array} \frac{-\pi^{2}}{16}\left\{\frac{(2)_{k}}{(1)_{k}}\right\}^{4}\right)\bigg]
=\frac{\pi^{2}}{60}\left(\frac{5}{\pi}-\frac{5}{\pi^{2}}-1\right),~~~~~~~~~~~~~~~~~~~~~
\end{align}
\begin{eqnarray}\label{N4}
\sum_{k=0}^{\infty}\bigg[\left\{\frac{(2)_{k}}{(1)_{k}}\right\}~{}_{0}F_{1}\left(\ \begin{array}{lll}\overline{~~~};\\ \frac{1}{2};\end{array} -\frac{\pi^{2}}{4}\left\{\frac{(2)_{k}}{(1)_{k}}\right\}^{4}\right)\bigg]
-2\sqrt{2}\sum_{k=0}^{\infty}\bigg[\left\{\frac{(2)_{k}}{(1)_{k}}\right\}^{2}{}_{1}F_{2}\left(\ \begin{array}{lll}1;\\
\frac{3}{4},\frac{5}{4};\end{array} \frac{-\pi^{2}}{4}\left\{\frac{(2)_{k}}{(1)_{k}}\right\}^{4} \right)\bigg]\nonumber+\\
+\pi\sum_{k=0}^{\infty}\bigg[\left\{\frac{(2)_{k}}{(1)_{k}}\right\}^{3}{}_{0}F_{1}\left(\ \begin{array}{lll}\overline{~~~~};\\
\frac{3}{2}~;\end{array} \frac{-\pi^{2}}{4}\left\{\frac{(2)_{k}}{(1)_{k}}\right\}^{4} \right)\bigg]
=\frac{\sqrt{2}-1}{4},~~~~~~~
\end{eqnarray}
\begin{eqnarray}\label{N5}
\sum_{k=0}^{\infty}\bigg[\left\{\frac{(2)_{k}}{(1)_{k}}\right\}~{}_{0}F_{1}\left(\ \begin{array}{lll}\overline{~~~};\\ \frac{1}{2};\end{array} \frac{-\pi^{2}}{16}\left\{\frac{(2)_{k}}{(1)_{k}}\right\}^{4}\right)\bigg]
-2\sum_{k=0}^{\infty}\bigg[\left\{\frac{(2)_{k}}{(1)_{k}}\right\}^{2}{}_{1}F_{2}\left(\ \begin{array}{lll}1;\\
\frac{3}{4},\frac{5}{4};\end{array} \frac{-\pi^{2}}{16}\left\{\frac{(2)_{k}}{(1)_{k}}\right\}^{4} \right)\bigg]\nonumber+\\
+\frac{\pi}{2}\sum_{k=0}^{\infty}\bigg[\left\{\frac{(2)_{k}}{(1)_{k}}\right\}^{3}{}_{0}F_{1}\left(\ \begin{array}{lll}\overline{~~~};\\
\frac{3}{2}~;\end{array} \frac{-\pi^{2}}{16}\left\{\frac{(2)_{k}}{(1)_{k}}\right\}^{4} \right)\bigg]=\frac{1}{4},~~~~~~~~~~~~
\end{eqnarray}
\begin{eqnarray}\label{N6}
\sum_{k=0}^{\infty}\bigg[\left\{\frac{(2)_{k}}{(1)_{k}}\right\}~{}_{0}F_{1}\left(\ \begin{array}{lll}\overline{~~~};\\ \frac{1}{2};\end{array} \frac{-\pi^{2}}{64}\left\{\frac{(2)_{k}}{(1)_{k}}\right\}^{4}\right)\bigg]
-\sqrt{2}\sum_{k=0}^{\infty}\bigg[\left\{\frac{(2)_{k}}{(1)_{k}}\right\}^{2}{}_{1}F_{2}\left(\ \begin{array}{lll}1;\\
\frac{3}{4},\frac{5}{4};\end{array} \frac{-\pi^{2}}{64}\left\{\frac{(2)_{k}}{(1)_{k}}\right\}^{4} \right)\bigg]\nonumber+\\
+\frac{\pi}{4}\sum_{k=0}^{\infty}\bigg[\left\{\frac{(2)_{k}}{(1)_{k}}\right\}^{3}{}_{0}F_{1}\left(\ \begin{array}{lll}\overline{~~~~};\\
\frac{3}{2}~;\end{array} \frac{-\pi^{2}}{64}\left\{\frac{(2)_{k}}{(1)_{k}}\right\}^{4} \right)\bigg]
=\frac{3\sqrt{2}-2}{4},~~~~~~~~~~~~~
\end{eqnarray}
\begin{eqnarray}\label{N7}
\sum_{k=0}^{\infty}\bigg[\left\{\frac{(2)_{k}}{(1)_{k}}\right\}~{}_{0}F_{1}\left(\ \begin{array}{lll}\overline{~~~};\\ \frac{1}{2};\end{array} \frac{-\pi^{2}}{144}\left\{\frac{(2)_{k}}{(1)_{k}}\right\}^{4}\right)\bigg]
-\frac{2\sqrt{3}}{3}\sum_{k=0}^{\infty}\bigg[\left\{\frac{(2)_{k}}{(1)_{k}}\right\}^{2}{}_{1}F_{2}\left(\ \begin{array}{lll}1;\\
\frac{3}{4},\frac{5}{4};\end{array} \frac{-\pi^{2}}{144}\left\{\frac{(2)_{k}}{(1)_{k}}\right\}^{4} \right)\bigg]\nonumber+\\
+\frac{\pi}{6}\sum_{k=0}^{\infty}\bigg[\left\{\frac{(2)_{k}}{(1)_{k}}\right\}^{3}{}_{0}F_{1}\left(\ \begin{array}{lll}\overline{~~~~};\\
\frac{3}{2}~;\end{array} \frac{-\pi^{2}}{144}\left\{\frac{(2)_{k}}{(1)_{k}}\right\}^{4} \right)\bigg]
=\frac{13\sqrt{3}-12}{12},~~~~~~~~~~~~~~~
\end{eqnarray}
\begin{eqnarray}\label{N8}
\sum_{k=0}^{\infty}\bigg[\left\{\frac{(2)_{k}}{(1)_{k}}\right\}~{}_{0}F_{1}\left(\ \begin{array}{lll}\overline{~~~};\\ \frac{1}{2};\end{array} -\pi^{2}\left\{\frac{(2)_{k}}{(1)_{k}}\right\}^{4}\right)\bigg]
-4\sum_{k=0}^{\infty}\bigg[\left\{\frac{(2)_{k}}{(1)_{k}}\right\}^{2}{}_{1}F_{2}\left(\ \begin{array}{lll}1;\\
\frac{3}{4},\frac{5}{4};\end{array} -\pi^{2}\left\{\frac{(2)_{k}}{(1)_{k}}\right\}^{4} \right)\bigg]\nonumber+\\
+2\pi\sum_{k=0}^{\infty}\bigg[\left\{\frac{(2)_{k}}{(1)_{k}}\right\}^{3}{}_{0}F_{1}\left(\ \begin{array}{lll}\overline{~~~~};\\
\frac{3}{2}~;\end{array} -\pi^{2}\left\{\frac{(2)_{k}}{(1)_{k}}\right\}^{4} \right)\bigg]=\frac{1}{8\pi},~~~~~~~~~~~~~~~~
\end{eqnarray}
\begin{align}\label{N9}
\sum_{k=0}^{\infty}\bigg[\left\{\frac{(2)_{k}}{(1)_{k}}\right\}~{}_{0}F_{1}\left(\ \begin{array}{lll}\overline{~~~};\\ \frac{1}{2};\end{array} \frac{-25\pi^{2}}{16}\left\{\frac{(2)_{k}}{(1)_{k}}\right\}^{4}\right)\bigg]
-2\sqrt{5}\sum_{k=0}^{\infty}\bigg[\left\{\frac{(2)_{k}}{(1)_{k}}\right\}^{2}{}_{1}F_{2}\left(\ \begin{array}{lll}1;\\
\frac{3}{4},\frac{5}{4};\end{array} \frac{-25\pi^{2}}{16}\left\{\frac{(2)_{k}}{(1)_{k}}\right\}^{4} \right)\bigg]\nonumber+\\
+\frac{5\pi}{2}\sum_{k=0}^{\infty}\bigg[\left\{\frac{(2)_{k}}{(1)_{k}}\right\}^{3}{}_{0}F_{1}\left(\ \begin{array}{lll}\overline{~~~~};\\
\frac{3}{2}~;\end{array} \frac{-25\pi^{2}}{16}\left\{\frac{(2)_{k}}{(1)_{k}}\right\}^{4} \right)\bigg]
=\frac{8\sqrt{5}-15}{100}.~~~~~~~~~~~~~~~~~~~~~~~~
\end{align}
The results (\ref{N1}) to (\ref{N3}) are obtained by putting $m=1, n=\frac{1}{2}$ ; $m=1, n=2$ and $m=2, n=2$ in the equations (\ref{K2}) and (\ref{F2}) and finally comparing with equations (\ref{s1}), (\ref{s2}) and (\ref{s3}). When $m=0$ with $n=1,~ 2, ~4, ~6, ~\frac{1}{2},~ \frac{2}{5}$ in the equations (\ref{K2}) and (\ref{F2}) and comparing with equations (\ref{P6}), (\ref{P7}), (\ref{P8}), (\ref{P9}), (\ref{P10}) and (\ref{P11}), we get the remaining results (\ref{N4}) to (\ref{N9}) respectively.
In view of the hypergeometric functions (\ref{E}), (\ref{E1}) and (\ref{E2}), we can express the above results (\ref{N4}) to (\ref{N9}) in terms of cosine, sine and Lommel functions.\\
Our results (\ref{N1}) to (\ref{N9}) are convergent in view of the convergence condition of ${}_{p}F_{q}(\cdot)$ series, when $p\leq q$, and for all $|z|<\infty$.
\section{Conclusion}
Here, we have described some infinite Fourier cosine transforms of Ramanujan. Thus certain Ramanujan's integrals, which may be different from those of presented here, can also be evaluated in a similar way. The results established above may be of significant in nature. We conclude our observation by remarking that various new results and applications can be obtained from our general theorem by appropriate choice of parameters $\upsilon,\lambda,b,c,y$ and bounded sequence $\{\Theta(k)\}_{k=0}^{\infty}$ in $\textbf{I}^{*}_{C}(\upsilon,b,c,\lambda,y)$. This work is in continuation to our earlier work \cite{Q8} on infinite Fourier sine transforms.


\section*{References}
\begin{thebibliography}{11}
\bibitem{A2} Agarwal, R. P. ; \textit{Resonance of Ramanujan's Mathematics}. Vol.I, New Age International (p) Limited Publisher, New Delhi, 1996.
\bibitem{A3} Agarwal, R. P. ; \textit{Resonance of Ramanujan's Mathematics}. Vol.II, New Age International (p) Limited Publisher, New Delhi, 1996.
\bibitem{A4} Agarwal, R. P. ; \textit{Resonance of Ramanujan's Mathematics}. Vol.III, New Age International (p) Limited Publisher, New Delhi, 1999.
\bibitem{A11} Andrews, G. E. and Berndt, B. C. ; \textit{Ramanujan's Lost Notebook.Part I}. Springer-Verlag, New York, 2005.
\bibitem{A12} Andrews, G. E. and Berndt, B. C. ; \textit{Ramanujan's Lost Notebook.Part II}. Springer-Verlag, New York, 2009.
\bibitem{A13} Andrews, G. E. and Berndt, B. C. ; \textit{Ramanujan's Lost Notebook.Part III}. Springer-Verlag, New York, 2012.
\bibitem{A14} Andrews, G. E. and Berndt, B. C. ; \textit{Ramanujan's Lost Notebook.Part IV}. Springer-Verlag, New York, 2013.
\bibitem{A15} Andrews, L.C. and Shivamoggi, B. K. ; \textit{Integral Transforms for Engineers}. Prentice-Hall of India, New Delhi, 2003.
\bibitem{B3} Beerends, R. J., Morsche, H. G. ter, Berg, J. C. Van den and Vri, E. M. Van de ; \textit{Fourier and Laplace Transforms}. Translated from Dutch by R.J. Beerends, Cambridge University Press, 2003.
\bibitem{B10} Berndt, B. C. ; \textit{Ramanujan's Notebooks.Part I}. Springer-Verlag, New York, 1985.
\bibitem{B11} Berndt, B. C. ; \textit{Ramanujan's Notebooks.Part II}. Springer-Verlag, New York, 1989.
\bibitem{B12} Berndt, B. C. ; \textit{Ramanujan's Notebooks.Part III}. Springer-Verlag, New York, 1991.
 \bibitem{B13} Berndt, B. C. ; \textit{Ramanujan's Notebooks.Part IV}. Springer-Verlag, New York, 1994.
 \bibitem{B14} Berndt, B. C. ; \textit{Ramanujan's Notebooks.Part V}. Springer-Verlag, New York, 1998.
\bibitem{C1} Campbell, G.A. and Foster, R.M. ; \textit{Fourier Integrals for Practical Applications.} Van Nostrand, New york, 1948.
\bibitem{C2} Carlson, B. C.; Some extensions of Lardner's relations between ${}_{0}F_{3}$ and Bessel functions. \textit{SIAM J. Math. Anal.} \textbf{1}(1970), 232-242.
\bibitem{C4} Carslaw, H. S. ; \textit{Introduction to the theory of Fourier's series and integrals}. Macmillan and co., limited st. Martin's street, London, 1921.
\bibitem{D} Ditkin, V.A. and Prudnikov, A.P. ; \textit{Integral transforms and operational calculus}. Pergamon Press, Oxford, London, Frankfurt, 1965.
\bibitem{E1} Erd\'{e}lyi, A., Magnus, W., Oberhettinger, F. and Tricomi, F. G. ; \textit{Higher Transcendental Functions}. Vol.1. McGraw-Hill, New york, Toronto and London, 1953.
\bibitem{E2} Erd\'{e}lyi, A., Magnus, W., Oberhettinger, F. and Tricomi, F. G. ; \textit{Tables of Integral Transforms}. Vol. 1.  McGraw-Hill, New york, Toronto and London, 1954.
\bibitem{K5} Kilbas,A.A. and Saigo, M. ; \textit{H-Transforms: Theory and Applications (Analytical Methods and Special Functions)}. CRC Press Company, Boca Raton, London, New York, Washington, D.C., Vol.\textbf{9}, 2004.
\bibitem{K6} Kilbas, A.A., Saigo, M. and Trujillo, J.J. ; On the generalized Wright function. \textit{Fract. Calc. Appl. Anal.} \textbf{5(4)}, 2002,  437-460.
\bibitem{L1} Lardner,T.J. ;  Relations between ${}_{0}F_{3}$ and Bessel functions. \textit{SIAM Rev}. \textbf{11}(1969) ,69-72.
\bibitem{M1} MacRobert, T. M. ; Integrals involving a modified Bessel function of the second kind and an E-function. \textit{Proc. Glasgow Math. Assoc.} \textbf{2}(1954), 93-96.
\bibitem{M2} Marichev, O.I. ; \textit{Handbook of Integral Transforms of Higher Transcendental Functions: Theory and algorithmic Tables}. Ellis. Horwood Ltd. John Wiley, 1983.
\bibitem{M} Meyer, J.L.; A generalization of an integral of Ramanujan. \textit{The Ramanujan Journal}, \textbf{14}(2007), 79-88.
\bibitem{O1} Oberhettinger, F. ; \textit{Tables of Mellin Transforms}. Springer-Verlag, Berlin, Heidelberg, New York, 1974.
\bibitem{o} Oberhettinger, F. ; \textit{Tables of Fourier Transforms and Fourier Transforms of Distributions}. Springer Verlag, Berlin, 1990.
\bibitem{P2} Prudnikov, A. P., Brychkov, Y.A. and Marichev, O.I. ;  \textit{Integrals and Series: Volume 1: Elementary Functions}. Taylor and Francis. 1986.
\bibitem{P3} Prudnikov, A. P., Brychkov, Y.A. and Marichev, O.I. ;  \textit{Integrals and Series: Volume 2: Special Functions}. Taylor and Francis, 1986.
\bibitem{P4} Prudnikov, A. P., Brychkov, Y.A. and Marichev, O.I. ;  \textit{Integrals and Series: Volume 3: More special functions}. Gordon and Breach Science Publishers, 1990.
\bibitem{P5} Prudnikov, A. P., Brychkov, Y.A. and Marichev, O.I. ;  \textit{Integrals and Series: Volume 4: Direct Laplace transforms}. Gordon and Breach Science Publishers, 1992.
\bibitem{P6} Prudnikov, A. P., Brychkov, Y.A. and Marichev, O.I. ;  \textit{Integrals and Series: Volume 5: Inverse Laplace transforms}. Gordon and Breach Science Publishers, 1992.
\bibitem{Q12} Qureshi, M. I.  and Khan, I. H. ; Ramanujan Integrals and other definite Integrals Associated with Gaussian hypergeometric functions. \textit{South East Asian Journal of Mathematics and Mathematical Sciences}, \textbf{4(1)}, 2005, 39-52.
\bibitem{Q9} Qureshi, M. I., Quraishi, K. A. and Pal, R. ; A Class of Hypergeometric Generalizations of an Integral of Srinivasa Ramanujan. \textit{Asian J. Current Engineering Math.} \textbf{2(3)}, 2013, 190-194.
\bibitem{Q4} Qureshi, M. I., Quraishi, K. A. and Pal, R. ;  Some Applications of Celebrated Master Theorem of Ramanujan. \textit{British Journal of Mathematics and Computer Science}, \textbf{4(20)}, 2014, 2862-2871.
\bibitem{Q2} Qureshi, M. I. and  Dar, S. A. ; Evaluation of Some Definite Integrals of Ramanujan, using Hypergeometric Approach. \textit{Palestine Journal of Mathematics}, \textbf{6(1)}, 2017, 1-3.
\bibitem{Q8} Qureshi, M. I. and Dar, S. A. ; Generalizations of Ramanujan's integral associated with infinite Fourier sine transforms in terms of hypergeometric functions and its applications. \textit{Communicated for publication}.
\bibitem{R1} Ramanujan, S. ; Some definite integrals connected with Gauss's sums. \textit{Mess. Math}., XLIV(1915), 75-86.
\bibitem{R4} Ramanujan, S. ; Some definite integrals. \textit{J. Indian Math. Soc.} \textbf{11}(1915), 81-87.
\bibitem{R0} Ramanujan's, S. ; \textit{Notebooks Vol. 1}. Tata Institute of Fundamental Research, Bombay, 1957.
\bibitem{R5} Ramanujan's, S. ; \textit{Notebooks Vol. 2}. Tata Institute of Fundamental Research, Bombay, 1957.
\bibitem{R3} Ramanujan, S. ;  \textit{The Lost Notebook and Other Unpublished Papers}. Narosa, New Delhi, 1988.
\bibitem{R2} Ramanujan, S. ; \textit{Collected papers}. Cambridge University Press, Cambridge, 1927, reprinted by Chelsea, New York, 1962, reprinted by Amer. Math. Soc., Providence, RI, 2000.
\bibitem{S6} Sharma, B. L. ;  A formula for hypergeometric series and its application. \textit{An. Univ. Timisoara Ser. Sti. Mat.} \textbf{12}(1974), 145-154.
\bibitem{S7} Sharma, B. L. ;  New generating functions for the Gegenbauer and the Hermite polynomial. \textit{Simon Stevin}, \textbf{54}(1980), 129-134.
\bibitem{s3} Sneddon, Ian, N. ; \textit{Fourier Transforms}. McGraw Hill Book Company, Inc, Newyork, 1951.
\bibitem{s4} Sneddon, Ian, N. ; \textit{The use of Integral Transforms}.  McGraw Hill Book Company, Inc, Newyork, 1972.
\bibitem{S5} Srivastava, H. M. ; A note on certain identities involving generalized hypergeometric series. \textit{Nederl. Akad. Wetensch. Proc. Ser. A 82=Indag. Math.} \textbf{41}(1979), 191-201.
\bibitem{s1} Srivastava, H. M. and Manocha, H. L. ; \textit{A Treatise on Generating Functions}. Halsted Press (Ellis Horwood Limited, Chichester, U.K.), John Wiley and Sons, New york, Chichester,Brisbane and Toronto, 1984.
\bibitem{s2} Srivastava, H. M., Qureshi, M.I., Singh, A. and Arora, A. ; A family of hypergeometric integrals associated with Ramanujan's integral formula. \textit{Advanced Studies in Contemporary Mathematics (Kyungshang)}, \textbf{18(2)}, 2009, 113-125.
\bibitem{T1} Titchmarsh, E. C. ; \textit{Introduction to the Theory of Fourier Integrals}. Clarendon Press, Second Edition , 1948.
\bibitem{W1} Wiener, N. ; \textit{The Fourier Integral and certain of its applications}. Dover publications, Newyork, 1951.
\bibitem{W2} Wright, E.M. ;  The asymptotic expansion of the generalized hypergeometric function. \textit{J. London Math. Soc.} \textbf{10}(1935), 286-293.
\bibitem{W4} Wright, E.M. ;  The asymptotic expansion of the generalized hypergeometric function. \textit{Proc. London Math. Soc.} \textbf{2(46)}, 1940, 389-408.
\end{thebibliography}
\end{document}